\documentclass[a4paper,12pt]{amsart}
\usepackage{amsfonts,amsmath,amsthm,amssymb}
\usepackage{enumitem}
\usepackage{mathrsfs}
\usepackage[margin=1in]{geometry}
\usepackage[colorlinks=true,linkcolor=blue,filecolor=blue,citecolor=red]{hyperref}

\usepackage{cleveref}

\theoremstyle{plain}
\newtheorem{theorem}{Theorem}[section]

\newtheorem{lemma}[theorem]{Lemma}

\theoremstyle{definition}

\begin{document}
\title[Modular equations involving Rogers-Ramanujan continued fraction]{A remark on modular equations involving Rogers-Ramanujan continued fraction via $5$-dissections}
\author[Russelle Guadalupe]{Russelle Guadalupe}
\address{Institute of Mathematics, University of the Philippines-Diliman\\
Quezon City 1101, Philippines}
\email{rguadalupe@math.upd.edu.ph}

\renewcommand{\thefootnote}{}



\renewcommand{\thefootnote}{\arabic{footnote}}
\setcounter{footnote}{0}

\begin{abstract}
In this paper, we study the $5$-dissections of certain Ramanujan's theta functions, particularly $\psi(q)\psi(q^2), \varphi(-q)$ and $\varphi(-q)\varphi(-q^2)$, and derive an identity for $q(q;q)_{\infty}^6/(q^5;q^5)_{\infty}^6$ in terms of certain products of the Rogers-Ramanujan continued fraction $R(q)$. Using this identity, we give another proof of the modular equation involving $R(q), R(q^2)$ and $R(q^4)$, which was recorded by Ramanujan in his lost notebook, and establish modular equations involving $R(q), R(q^2), R(q^4), R(q^8)$ and $R(q^{16})$.
\end{abstract}

\maketitle

\section{Introduction}\label{sec1}

Throughout this paper, we define $(a;q)_{\infty} := \prod_{n=0}^\infty (1-aq^n)$ for complex numbers $a$ and $q$ with $|q| < 1$ and
\begin{align*}
(a_1,a_2,\ldots,a_m;q)_\infty:= (a_1;q)_\infty (a_2;q)_\infty\cdots (a_m;q)_\infty. 
\end{align*}

In 1894, Rogers \cite{rogers} first studied certain properties of the continued fraction 
\begin{align*}
R(q) := \dfrac{1}{1+\dfrac{q}{1+\dfrac{q^2}{1+\dfrac{q^3}{1+\cdots}}}},
\end{align*}
now known as the Rogers-Ramanujan continued fraction, and gave its infinite product expansion
\begin{align*}
R(q) = \dfrac{(q,q^4;q^5)_\infty}{(q^2,q^3;q^5)_\infty}.
\end{align*}

Almost two decades later, Ramanujan rediscovered $R(q)$ in his first letter \cite[pp. 21--30]{ramlet} to Hardy, where he claimed that if $u=q^{1/5}R(q)$ and $v = qR(q^5)$, then 
\begin{align}
v^5 = u\cdot \dfrac{1-2u+4u^2-3u^3+u^4}{1+3u+4u^2+2u^3+u^4}.\label{eq11}
\end{align}

In his notebooks \cite{ramnot,ramlos}, Ramanujan recorded several remarkable identities involving $R(q)$, including the following modular equations

\begin{align}
\dfrac{R_2-R_1^2}{R_2+R_1^2}&=qR_1R_2^2,\label{eq12}\\
\dfrac{R_4-R_1^2R_2}{R_4+R_2^2}&=qR_1R_4,\label{eq13}
\end{align}
where $R_m:=R(q^m)$. The identity (\ref{eq12}) appeared in Ramanujan's notebook \cite[p. 326]{ramnot}, and its proof, together with that of (\ref{eq11}), was later provided by Rogers \cite{rogers2} using the theory of elliptic functions. Other proofs of (\ref{eq12}) were given by Berndt \cite[p. 12, Entry 1]{berndt5} using a modular equation of degree five, and Hirschhorn \cite[p. 367]{hirsc} using two $q$-series identities involving $R_1$ and $R_2$. On the other hand, the identity (\ref{eq13}) appeared in Ramanujan's lost notebook \cite[p. 205]{ramlos}, and was later established by Andrews and Berndt \cite[pp. 24-25]{andber1} using (\ref{eq12}).

Recently, Tang \cite{tang} gave another proofs of (\ref{eq12}) and (\ref{eq13}) using the $5$-dissections of the Ramanujan's theta functions $(q;q)_\infty, (q;q)_\infty^{-1}$ and $\psi(q)$. Tang also gave several $q$-series identities involving theta functions, which includes instances of $(q;q)_{\infty}^6/(q^5;q^5)_{\infty}^6$, as byproducts. 

Motivated by the work of Tang \cite{tang}, we explore in this paper the $5$-dissections of the Ramanujan's theta functions $\psi(q)\psi(q^2), \varphi(-q)$ and $\varphi(-q)\varphi(-q^2)$, and derive an identity for $q(q;q)_{\infty}^6/(q^5;q^5)_{\infty}^6$ in terms of $R_1, R_2$ and $R_4$, as shown in our main result below.

\begin{theorem}\label{thm11}
We have
\begin{align}
	q\dfrac{(q;q)_{\infty}^6}{(q^5;q^5)_{\infty}^6}=\dfrac{P(R_1,R_2,R_4)}{R_1^4R_2R_4^2},\label{eq14}
\end{align}
where
\begin{align*}
	P(a,b,c) &:= q^4(-2 a^8 b^2 c^3 + a^8 c^4 - 2 a^6 b c^4 - 5 a^4 b^2 c^4) + q^3(-a^7 b^2 c^2 - 2 a^7 c^3+ 6 a^5 b c^3\\
	&- 6 a^3 b^2 c^3 + 3 a^3 c^4 - a b c^4) + q^2(2 a^8 b c + 4 a^6 b^2 c - 2 a^6 c^2 + 5 a^4 b c^2+ 2 a^2 b^2 c^2\\
	&+ 4 a^2 c^3 - 2 b c^3)+q(a^7 b + 3 a^5 b^2 + 6 a^5 c + 6 a^3 b c + 2 a b^2 c - a c^2)+5 a^4 - 2 a^2 b\\
	&- b^2 - 2 c.
\end{align*}
\end{theorem}

Using Theorem \ref{thm11}, we provide another proof of the identity (\ref{eq13}), and establish the following modular equations involving $R(q), R(q^2), R(q^4), R(q^8)$ and $R(q^{16})$, in the spirit of (\ref{eq12}) and (\ref{eq13}).

\begin{theorem}\label{thm12}
We have
\begin{align}
	\dfrac{(R_2-R_1^2)(R_8-R_1^2R_2R_4)}{(R_2-R_1^2)(R_8+R_4^2)+R_2^3(R_2+R_1^2)^2} &= q^2R_1^2R_8,\label{eq15}\\
	\dfrac{(R_2-R_1^2)^3(R_4-R_2^2)(R_{16}-R_1^2R_2R_4R_8)}{(R_2-R_1^2)^3(R_4-R_2^2)(R_{16}+R_8^2)+Q(R_1,R_2,R_4)} &=q^4R_1^2R_2R_{16},\label{eq16}
\end{align}
where 
\begin{align*}
	Q(a,b,c):=a^2b^6(b+a^2)^4(c-b^2)+c^3(c+b^2)^2(b-a^2)^3.
\end{align*}
\end{theorem}

The rest of the paper is organized as follows. In Section \ref{sec2}, we recall the $5$-dissections of $(q;q)_\infty, (q;q)_\infty^{-1}$ and $\psi(q)$, and deduce the $5$-dissection of $\psi(q)\psi(q^2)$. We then utilize these dissections to prove Theorem \ref{thm11}, and establish (\ref{eq13}) via an identity of Ramanujan. In Sections \ref{sec3} and \ref{sec4}, we apply Theorem \ref{thm11} and the $5$-dissections of the theta functions $\varphi(-q)$ and $\varphi(-q)\varphi(-q^2)$ to prove Theorem \ref{thm12}, and give a short remark on the proof of this theorem. We perform most of our computations via \textit{Mathematica}.

\section{Proofs of Theorem \ref{thm11} and identity (\ref{eq13})}\label{sec2}

We recall the Ramanujan's general theta function \cite[p. 34, (18.1)]{berndt3}
\begin{align*}
f(a,b) = \sum_{n=-\infty}^{\infty} a^{n(n+1)/2}b^{n(n-1)/2} = (-a,-b,ab;ab)_\infty
\end{align*}
valid for complex numbers $a$ and $b$ with $|ab| < 1$, where the last equality follows from the Jacobi triple product identity \cite[p. 35, Entry 19]{berndt3}. We need three special cases of $f(a,b)$, namely \cite[p. 36, Entry 22]{berndt3}
\begin{align*}
\varphi(q) &:= f(q,q) = \sum_{n=-\infty}^\infty q^{n^2} = (-q;q^2)_\infty^2(q^2;q^2)_\infty,\\
\psi(q) &:= f(q,q^3) =\sum_{n=0}^\infty q^{n(n+1)/2} = \dfrac{(q^2;q^2)_\infty^2}{(q;q)_\infty},\\
f(-q) &:= f(-q,-q^2) = \sum_{n=-\infty}^\infty (-1)^nq^{n(3n-1)/2} = (q;q)_\infty.
\end{align*}

We next introduce the $5$-dissections of $(q;q)_\infty, (q;q)_\infty^{-1}$ and $\psi(q)$, as shown in the following lemma.

\begin{lemma}\label{lem21}
We have
\begin{align}
	(q;q)_\infty &= (q^{25};q^{25})_\infty\left(\dfrac{1}{R_5}-q-q^2R_5\right),\label{eq21}\\
	\dfrac{1}{(q;q)_\infty} &= \dfrac{(q^{25};q^{25})_\infty^5}{(q^5;q^5)_\infty^6}\left(\dfrac{1}{R_5^4}+\dfrac{q}{R_5^3}+\dfrac{2q^2}{R_5^2}+\dfrac{3q^3}{R_5}+5q^4-3q^5R_5+2q^6R_5^2-q^7R_5^3+q^8R_5^4\right),\label{eq22}\\
	\psi(q) &= (-q^{10},-q^{15},q^{25};q^{25})_\infty+q(-q^5,-q^{20},q^{25};q^{25})_\infty +q^3\psi(q^{25}).\label{eq23}
\end{align}
\end{lemma}

\begin{proof}
See \cite[pp. 86-88]{hirsc} for the proof of (\ref{eq21}), and \cite[p. 184]{andber3} for the proof of (\ref{eq22}). On the other hand, identity (\ref{eq23}) follows from \cite[p. 49, Cor. (ii)]{berndt3}.
\end{proof}

We now in a position to prove Theorem \ref{thm11}. Consequently, we offer another proof of the identity (\ref{eq13}) using Theorem \ref{thm11}.

\begin{proof}[Proof of Theorem \ref{thm11}]
By (\ref{eq23}), we have the $5$-dissection of $\psi(q)\psi(q^2)$ given by 
\begin{align}
	\psi(q)\psi(q^2) &= \left((-q^{10},-q^{15},q^{25};q^{25})_\infty+q(-q^5,-q^{20},q^{25};q^{25})_\infty +q^3\psi(q^{25})\right)\nonumber\\
	&\times \left((-q^{20},-q^{30},q^{50};q^{50})_\infty+q^2(-q^{10},-q^{40},q^{50};q^{50})_\infty +q^6\psi(q^{50})\right)\nonumber\\
	&=: A_0(q^5)+qA_1(q^5)+q^2A_2(q^5)+q^3A_3(q^5)+q^9\psi(q^{25})\psi(q^{50})\label{eq24}
\end{align}
for some $A_k(q)\in\mathbb{Z}[[q]]$ for $k\in\{0,\ldots 3\}$. On the other hand, we apply (\ref{eq21}) and (\ref{eq22}) on $\psi(q)\psi(q^2)$ and obtain
\begin{align}
	\psi(q)\psi(q^2) &= \dfrac{(q^2;q^2)_\infty(q^4;q^4)_\infty^2}{(q;q)_\infty}\nonumber\\
	&=(q^{50};q^{50})_\infty(q^{100};q^{100})_\infty^2\left(\dfrac{1}{R_{10}}-q^2-q^4R_{10}\right)\left(\dfrac{1}{R_{20}}-q^4-q^8R_{20}\right)^2\nonumber\\
	&\times \dfrac{(q^{25};q^{25})_\infty^5}{(q^5;q^5)_\infty^6}\left(\dfrac{1}{R_5^4}+\dfrac{q}{R_5^3}+\dfrac{2q^2}{R_5^2}+\dfrac{3q^3}{R_5}+5q^4-3q^5R_5+2q^6R_5^2-q^7R_5^3+q^8R_5^4\right).\label{eq25}
\end{align}
Comparing the terms of (\ref{eq24}) and (\ref{eq25}) involving $q^{5n+4}$, dividing both sides by $q^4$, and replacing $q^5$ with $q$, we get
\begin{align*}
	q\psi(q^5)\psi(q^{10}) = (q^{10};q^{10})_\infty(q^{20};q^{20})_\infty^2 \cdot \dfrac{(q^5;q^5)_\infty^5}{(q;q)_\infty^6}\dfrac{P(R_1,R_2,R_4)}{R_1^4R_2R_4^2},
\end{align*}
which yields (\ref{eq14}).
\end{proof}

\begin{proof}[Proof of (\ref{eq13})]
We require the following identity of Ramanujan \cite[p. 213, (16)]{ramp}, \cite[p. 267, (11.6)]{berndt3} given by
\begin{align}
	\dfrac{1}{qR_1^5}-11-qR_1^5 = \dfrac{(q;q)_\infty^6}{q(q^5;q^5)_\infty^6}.\label{eq26}
\end{align}
Comparing equation (\ref{eq14}) of Theorem \ref{thm11} and (\ref{eq26}) yields
\begin{align*}
	\dfrac{P(R_1,R_2,R_4)}{R_1^4R_2R_4^2} = q^2\left(\dfrac{1}{qR_1^5}-11-qR_1^5\right),
\end{align*}
which reduces to 
\begin{align}
	\dfrac{(1+qR_1R_4)^2X(R_1,R_2,R_4)}{R_1^5R_2R_4^2}=0, \label{eq27}
\end{align}
where
\begin{align}
	X(a,b,c) &:= q^2(-2 a^7 b^2 c + a^7 c^2 - 2 a^5 b c^2 - 5 a^3 b^2 c^2)+q(a^8 b + 3 a^6 b^2 - 4 a^6 c + 10 a^4 b c\nonumber\\
	&+ 4 a^2 b^2 c + 3 a^2 c^2 - b c^2)+5 a^5 - 2 a^3 b - a b^2 - 2 a c.\label{eq28}
\end{align}
We now write
\begin{align}
	R_1X(R_1,R_2,R_4) &= (-R_1^2R_2+R_4-qR_1R_2^2R_4-qR_1R_4^2)Y(R_1,R_2)\nonumber\\
	&-(R_1^2-R_2+qR_1^3R_2^2+qR_1R_2^3)(-5R_1^4-R_4+3qR_1^5R_4-5qR_1^3R_2R_4),\label{eq29}
\end{align}
where $Y(a,b) = q(-a^7 + 2 a^5 b + 5 a^3 b^2)-3 a^2 + b$. We note that (\ref{eq12}) is equivalent to 
\begin{align}
	R_1^2-R_2+qR_1^3R_2^2+qR_1R_2^3 =0. \label{eq210}
\end{align}
Because 
\begin{align*}
	1+qR_1R_4 &= 1 + q - q^2 + q^3 - 2 q^5 + O(q^6),\\
	Y(R_1,R_2) &= -2 + 12 q - 28 q^2 + 26 q^3 + 26 q^4 - 144 q^5+O(q^6)
\end{align*}
are not identically zero, it follows from (\ref{eq27}), (\ref{eq29}) and (\ref{eq210}) that
\begin{align}
	-R_1^2R_2+R_4-qR_1R_2^2R_4-qR_1R_4^2=0,\label{eq211}
\end{align}
which transforms into (\ref{eq13}) after rearranging terms.
\end{proof}

\section{Proof of Theorem \ref{thm12}}\label{sec3}

To prove Theorem \ref{thm12}, we first record the following $5$-dissection of $\varphi(-q)$.

\begin{lemma}\label{lem31}
We have the $5$-dissection 
\begin{align}
	\varphi(-q) = \varphi(-q^{25})-2q(q^{15},q^{35},q^{50};q^{50})_\infty+2q^4(q^5,q^{45},q^{50};q^{50})_\infty.\label{eq31}
\end{align}
\end{lemma}

\begin{proof}
This follows from replacing $q$ with $-q$ in the second identity of \cite[p. 49, Cor. (i)]{berndt3}.
\end{proof}

\begin{proof}[Proof of Theorem \ref{thm12}]
We first prove (\ref{eq15}). We apply (\ref{eq21}) and (\ref{eq22}) on $\varphi(-q)$ and get
\begin{align}
	\varphi(-q) &= \dfrac{(q;q)_\infty^2}{(q^2;q^2)_\infty}\nonumber\\
	&=(q^{25};q^{25})_\infty^2\dfrac{(q^{50};q^{50})_\infty^5}{(q^{10};q^{10})_\infty^6}\left(\dfrac{1}{R_5}-q-q^2R_5\right)^2\nonumber\\
	&\times\left(\dfrac{1}{R_{10}^4}+\dfrac{q^2}{R_{10}^3}+\dfrac{2q^4}{R_{10}^2}+\dfrac{3q^6}{R_{10}}+5q^8-3q^{10}R_{10}+2q^{12}R_{10}^2-q^{14}R_{10}^3+q^{16}R_{10}^4\right).\label{eq32}
\end{align}
Taking the terms of (\ref{eq31}) and (\ref{eq32}) involving $q^{5n}$ and replacing $q^5$ with $q$, we deduce that
\begin{align}
	\dfrac{(q^2;q^2)_\infty^6}{(q^{10};q^{10})_\infty^6}=\dfrac{A(R_1,R_2)}{R_1^2R_2^4},\label{eq33}
\end{align}
where 
\begin{align*}
	A(a,b) &= q^4a^4b^8+q^3(4 a^3 b^6 + 2 a b^7)+q^2(3 a^4 b^3 - 5 a^2 b^4 - 3 b^5)+q(2 a^3 b - 4 a b^2)+1.
\end{align*}
We now replace $q$ with $q^2$ in the coefficients of the polynomial $P(a,b,c)$ in Theorem \ref{thm11} to get a new polynomial $P_0(a,b,c)$ such that
\begin{align}
	q^2\dfrac{(q^2;q^2)_\infty^6}{(q^{10};q^{10})_\infty^6} = \dfrac{P_0(R_2,R_4,R_8)}{R_2^4R_4R_8^2}.\label{eq34}
\end{align}
Comparing (\ref{eq33}) and (\ref{eq34}) yields
\begin{align*}
	q^2\dfrac{A(R_1,R_2)}{R_1^2R_2^4}=\dfrac{P_0(R_2,R_4,R_8)}{R_2^4R_4R_8^2},
\end{align*}
which expands to 
\begin{align}
	\dfrac{B(R_1,R_2,R_4,R_8)}{R_1^2R_2^4R_4R_8^2}=0,\label{eq35}
\end{align}
where 
\begin{align}
	B(a,b,c,d) &:= q^8(2 a^2 b^8 c^2 d^3 - a^2 b^8 d^4 + 2 a^2 b^6 c d^4 + 5 a^2 b^4 c^2 d^4) + q^6(a^4 b^8 c d^2 + a^2 b^7 c^2 d^2\nonumber\\
	&+ 2 a^2 b^7 d^3 - 6 a^2 b^5 c d^3 + 6 a^2 b^3 c^2 d^3 - 3 a^2 b^3 d^4 + a^2 b c d^4) + q^5(4 a^3 b^6 c d^2\nonumber\\
	&+ 2 a b^7 c d^2) + q^4(-2 a^2 b^8 c d - 4 a^2 b^6 c^2 d + 2 a^2 b^6 d^2 + 3 a^4 b^3 c d^2 - 10 a^2 b^4 c d^2\nonumber\\
	&- 3 b^5 c d^2 - 2 a^2 b^2 c^2 d^2 - 4 a^2 b^2 d^3 + 2 a^2 c d^3)+q^3(2 a^3 b c d^2 - 4 a b^2 c d^2)\nonumber\\
	&+q^2(-a^2 b^7 c - 3 a^2 b^5 c^2 - 6 a^2 b^5 d - 6 a^2 b^3 c d - 2 a^2 b c^2 d + a^2 b d^2 + c d^2)-5 a^2 b^4\nonumber\\
	&+ 2 a^2 b^2 c + a^2 c^2 + 2 a^2 d.\label{eq36}
\end{align}
Replacing $q$ with $q^2$ in the coefficients of the polynomial (\ref{eq28}), we have a new polynomial $X_0(a,b,c)$ with $X_0(R_2,R_4,R_8)=0$ in view of (\ref{eq27}). We now express
\begin{align}
	R_2B(R_1,R_2,R_4,R_8) &= q^3R_2^2R_4R_8^2(R_1^2-R_2+qR_1^3R_2^2+qR_1R_2^3)Z_0(R_1,R_2)\nonumber\\
	&-X_0(R_2,R_4,R_8)(R_1^2 + 2 q^2 R_1^2 R_2 R_8 + q^4 R_1^2 R_2^2 R_8^2)\nonumber\\
	&+q^2R_1^2R_4R_8(R_2^2R_4-R_8+q^2R_2R_4^2R_8+q^2R_2R_8^2)\nonumber\\
	&+q^2R_2R_4R_8(-R_1^2R_2R_4+R_8-qR_1^3R_2R_8-qR_1R_2^2R_8\nonumber\\
	&-q^2R_1^2R_4^2R_8-q^2R_1^2R_8^2),\label{eq37}
\end{align}
where $Z_0(a,b)=q^2ab^5+3qb^3+3a$. Replacing $q$ with $q^2$ in (\ref{eq210}) gives
\begin{align}
	R_2^2R_4-R_8+q^2R_2R_4^2R_8+q^2R_2R_8^2=0.\label{eq38}
\end{align}
Since $q^2R_2R_4R_8$ is not identically zero and $X_0(R_2,R_4,R_8)=0$, we infer from (\ref{eq210}), (\ref{eq35}), (\ref{eq37}) and (\ref{eq38}) that 
\begin{align}
	-R_1^2R_2R_4+R_8-qR_1^3R_2R_8-qR_1R_2^2R_8-q^2R_1^2R_4^2R_8-q^2R_1^2R_8^2=0,\label{eq39}
\end{align}
which can be written as 
\begin{align}
	\dfrac{R_8-R_1^2R_2R_4}{R_4^2+R_8+
		\dfrac{R_2(R_1^2+R_2)}{qR_1}} = q^2R_1^2R_8.\label{eq310}
\end{align}
We know from (\ref{eq12}) that
\begin{align}
	\dfrac{1}{qR_1} = \dfrac{R_2^2(R_1^2+R_2)}{R_2-R_1^2},\label{eq311}
\end{align}
so plugging this equation into (\ref{eq310}) yields (\ref{eq15}). We next prove (\ref{eq16}); by (\ref{eq31}) we derive the $5$-dissection of $\varphi(-q)\varphi(-q^2)$ given by 
\begin{align}
	\varphi(-q)\varphi(-q^2) &= \left(\varphi(-q^{25})-2q(q^{15},q^{35},q^{50};q^{50})_\infty+2q^4(q^5,q^{45},q^{50};q^{50})_\infty\right)\nonumber\\
	&\times \left(\varphi(-q^{50})-2q^2(q^{30},q^{70},q^{100};q^{100})_\infty+2q^8(q^{10},q^{90},q^{100};q^{100})_\infty\right)\nonumber\\
	&=: \varphi(-q^{25})\varphi(-q^{50})+qC_1(q^5)+q^2C_2(q^5)+q^3C_3(q^5)+q^4C_4(q^5)\label{eq312}
\end{align}
for some $C_k(q)\in\mathbb{Z}[[q]]$ for $k\in \{1,\ldots,4\}$. We next use (\ref{eq21}) and (\ref{eq22}) on $\varphi(-q)\varphi(-q^2)$ and get
\begin{align}
	\varphi(-q)\varphi(-q^2) &= \dfrac{(q;q)_\infty^2(q^2;q^2)_\infty}{(q^4;q^4)_\infty}\nonumber\\
	&=(q^{25};q^{25})_\infty^2(q^{50};q^{50})_\infty\dfrac{(q^{100};q^{100})_\infty^5}{(q^{20};q^{20})_\infty^6}\left(\dfrac{1}{R_5}-q-q^2R_5\right)^2\left(\dfrac{1}{R_{10}}-q^2-q^4R_{10}\right)\nonumber\\
	&\times\left(\dfrac{1}{R_{20}^4}+\dfrac{q^4}{R_{20}^3}+\dfrac{2q^8}{R_{20}^2}+\dfrac{3q^{12}}{R_{20}}+5q^{16}-3q^{20}R_{20}+2q^{24}R_{20}^2-q^{28}R_{20}^3+q^{32}R_{20}^4\right).\label{eq313}
\end{align}
Extracting the terms of (\ref{eq312}) and (\ref{eq313}) involving $q^{5n}$ and replacing $q^5$ with $q$, we arrive at
\begin{align}
	\dfrac{(q^4;q^4)_\infty^6}{(q^{20};q^{20})_\infty^6} = \dfrac{F(R_1,R_2,R_4)}{R_1^2R_2R_4^4},\label{eq314}
\end{align}
where 
\begin{align*}
	F(a,b,c) &= -q^8a^4 b^2 c^8+q^7(2 a^3 b^2 c^7 + 2 a^3 c^8 + 2 a b c^8)+q^6(-2 a^4 b c^6 + 2 a^2 b^2 c^6 + a^2 c^7 + b c^7)\\
	&+q^5(6 a^3 b c^5 - 6 a b^2 c^5 - 4 a c^6)+q^4(-3 a^4 b^2 c^3 + 5 a^4 c^4 + 5 a^2 b c^4 - 5 b^2 c^4 - 3 c^5)\\
	&+q^3(-4 a^3 b^2 c^2 + 6 a^3 c^3 + 6 a b c^3)+q^2(-a^4 b c + a^2 b^2 c - 2 a^2 c^2 - 2 b c^2)\\
	&+q(-2 a^3 b + 2 a b^2 - 2 a c)+1.
\end{align*}
We replace $q$ with $q^4$ in the coefficients of the polynomial $P(a,b,c)$ in Theorem \ref{thm11}, yielding a new polynomial $P_1(a,b,c)$ such that
\begin{align}
	q^4\dfrac{(q^4;q^4)_\infty^6}{(q^{20},q^{20})_\infty^6}=\dfrac{P_1(R_4,R_8,R_{16})}{R_4^4R_8R_{16}^2}.\label{eq315}
\end{align}
We see from (\ref{eq314}) and (\ref{eq315}) that
\begin{align*}
	q^4\dfrac{F(R_1,R_2,R_4)}{R_1^2R_2R_4^4}=\dfrac{P_1(R_4,R_8,R_{16})}{R_4^4R_8R_{16}^2},
\end{align*}
which reduces to 
\begin{align}
	\dfrac{C(R_1,R_2,R_4,R_8,R_{16})}{R_1^2R_2R_4^4R_8R_{16}^2}=0,\label{eq316}
\end{align}
where
\begin{align*}
	C(a,b,c,d,e) &:= q^{16}(2 a^2 b c^8 d^2 e^3 - a^2 b c^8 e^4 + 2 a^2 b c^6 d e^4 + 5 a^2 b c^4 d^2 e^4)+q^{12}(-a^4 b^2 c^8 d e^2\\
	&+ a^2 b c^7 d^2 e^2 + 2 a^2 b c^7 e^3 - 6 a^2 b c^5 d e^3 + 6 a^2 b c^3 d^2 e^3 - 3 a^2 b c^3 e^4 + a^2 b c d e^4)\\
	&+q^{11}(2 a^3 b^2 c^7 d e^2 + 2 a^3 c^8 d e^2 + 2 a b c^8 d e^2)+q^{10}(-2 a^4 b c^6 d e^2 + 2 a^2 b^2 c^6 d e^2\\
	&+ a^2 c^7 d e^2 + b c^7 d e^2)+q^9(6 a^3 b c^5 d e^2 - 6 a b^2 c^5 d e^2 - 4 a c^6 d e^2)+q^8(-2 a^2 b c^8 d e\\
	&- 4 a^2 b c^6 d^2 e + 2 a^2 b c^6 e^2 - 3 a^4 b^2 c^3 d e^2 + 5 a^4 c^4 d e^2 - 5 b^2 c^4 d e^2 - 3 c^5 d e^2\\
	&- 2 a^2 b c^2 d^2 e^2 - 4 a^2 b c^2 e^3 + 2 a^2 b d e^3)+q^7(-4 a^3 b^2 c^2 d e^2 + 6 a^3 c^3 d e^2\\
	&+ 6 a b c^3 d e^2)+q^6(-a^4 b c d e^2 + a^2 b^2 c d e^2 - 2 a^2 c^2 d e^2 - 2 b c^2 d e^2)+q^5(-2 a^3 b d e^2\\ 
	&+ 2 a b^2 d e^2 - 2 a c d e^2)+q^4(-a^2 b c^7 d - 3 a^2 b c^5 d^2 - 6 a^2 b c^5 e - 6 a^2 b c^3 d e\\
	&-2 a^2 b c d^2 e + a^2 b c e^2 + d e^2)-5 a^2 b c^4 + 2 a^2 b c^2 d + a^2 b d^2 + 2 a^2 b e.
\end{align*}

We replace $q$ with $q^2$ in the coefficients of the polynomial (\ref{eq36}) to form a new polynomial $B_0(a,b,c,d)$, so that from (\ref{eq35}) we have $B_0(R_2,R_4,R_8,R_{16})=0$. We expand
\begin{align}
	R_2^3C(R_1,R_2,R_4,R_8,R_{16}) &= -q^4R_2R_8R_{16}^2(R_1^2-R_2+qR_1^3R_2^2+qR_1R_2^3)G_0(R_1,R_2,R_4)\nonumber\\
	&+R_1^2R_2^2B_0(R_2,R_4,R_8,R_{16}) +q^4 R_1^2R_4R_8R_{16} (1 + 2 q^4 R_4^5)\nonumber\\
	&(-R_2^2R_4R_8 + R_{16} - q^2 R_2^3 R_4 R_{16} - q^2 R_2 R_4^2 R_{16} - q^4 R_2^2 R_8^2 R_{16}\nonumber\\
	&- q^4 R_2^2 R_{16}^2)-q^4R_8R_{16}^2(R_2^2-R_4+q^2R_2^3R_4^2+q^2R_2R_4^3)\nonumber\\
	&G_1(R_1,R_2,R_4)+5q^6R_4^3R_8R_{16}^2(1+qR_1R_4)(-R_1^2R_2+R_4\nonumber\\
	&-qR_1R_2^2R_4-qR_1R_4^2)-q^4R_2R_4R_8R_{16}(1+2q^4R_4^5)\nonumber\\
	&(-R_1^2 R_2 R_4 R_8 + R_{16} - q R_1^3 R_2 R_{16}- q R_1 R_2^2 R_{16} - q^2 R_1^2 R_2^2 R_4 R_{16}\nonumber\\
	&- q^2 R_1^2 R_2^2 R_{16} - q^4 R_1^2 R_2 R_8^2 R_{16}- q^4 R_1^2 R_2 R_{16}^2),\label{eq317}
\end{align}
where
\begin{align*}
	G_0(a,b,c) &:= q^7a b^2 c^8+q^6(-2 b^2 c^7 - 3 c^8)+2q^5 a b c^6-6q^4b c^5+q^3(3 a b^2 c^3 - 5 a c^4)\\
	&+q^2(4 b^2 c^2 - 9 c^3)+qa b c+2b,\\
	G_1(a,b,c) &:= 2q^5 a b^3 c^5+q^4(3 a^2 c^5 - 3 b c^5)+10q^3a b^2 c^3+q^2(-5 a^2 b c^2 + 5 b^2 c^2 - 5 c^3)\\
	&-4qa b^3-a^2 + b.
\end{align*}
We also replace $q$ with $q^2$ in (\ref{eq210}) and (\ref{eq39}) so that
\begin{align}
	R_2^2-R_4+q^2R_2^3R_4^2+q^2R_2R_4^3 = 0, \label{eq318}
\end{align}
and 
\begin{align}
	-R_2^2R_4R_8 + R_{16} - q^2 R_2^3 R_4 R_{16} - q^2 R_2 R_4^2 R_{16} - q^4 R_2^2 R_8^2 R_{16}
	- q^4 R_2^2 R_{16}^2 =0. \label{eq319}
\end{align}

Since $q^4R_2R_4R_8R_{16}(1+2q^4R_4^5)$ is not identically zero and $B_0(R_2,R_4,R_8,R_{16})=0$, we deduce from (\ref{eq210}), (\ref{eq211}), (\ref{eq316}), (\ref{eq317}), (\ref{eq318}) and (\ref{eq319}) that
\begin{align*}
	-R_1^2 R_2 R_4 R_8 + R_{16} - q R_1^3 R_2 R_{16} &- q R_1 R_2^2 R_{16} - q^2 R_1^2 R_2^2 R_4 R_{16}\\
	&- q^2 R_1^2 R_2^2 R_{16} - q^4 R_1^2 R_2 R_8^2 R_{16}- q^4 R_1^2 R_2 R_{16}^2=0,
\end{align*}
which can be written as
\begin{align}
	\dfrac{R_{16}-R_1^2R_2R_4R_8}{R_{16}+R_8^2+\dfrac{R_1^2(R_1^2+R_2)}{q^3R_1^3}+\dfrac{R_4(R_2^2+R_4)}{q^2R_2}}=q^4R_1^2R_2R_{16}.\label{eq320}
\end{align}
We note that (\ref{eq318}) is equivalent to
\begin{align*}
	\dfrac{1}{q^2R_2}=\dfrac{R_4^2(R_2^2+R_4)}{R_4-R_2^2},
\end{align*}
so plugging this identity together with (\ref{eq311}) in (\ref{eq320}) and simplifying, we finally arrive at (\ref{eq16}) as desired.
\end{proof}

\section{Concluding remarks}\label{sec4}

We have established the modular equations involving $R(q), R(q^2), R(q^4), R(q^8)$ and $R(q^{16})$ using the $5$-dissections of $\psi(q)\psi(q^2), \varphi(-q)$ and $\varphi(-q)\varphi(-q^2)$, as illustrated by (\ref{eq15}) and (\ref{eq16}) of Theorem \ref{thm12}, by expressing $q(q;q)_\infty^6/(q^5;q^5)_\infty^6$ in terms of $R(q), R(q^2)$ and $R(q^4)$ as shown in Theorem \ref{thm11}. We would like to remark, however, that we also use (\ref{eq12}) and (\ref{eq13}) to arrive at (\ref{eq15}) and (\ref{eq16}). Hence, it would be interesting to find a proof of Theorem \ref{thm12} without utilizing  (\ref{eq12}) and (\ref{eq13}).

%

\bibliographystyle{amsplain}

\end{document}